\documentclass{amsart}
\newtheorem{theorem}{Theorem}

\newtheorem{corollary}{Corollary}
\theoremstyle{definition}

\theoremstyle{remark}

\newcommand{\Z}{{\bf Z}}
\newcommand{\ra}{\rightarrow}
\newcommand{\fq}{\mathbf{F}_{q}}
\newcommand{\norm}{\mathrm{N}}

\begin{document}
\title{The Polynomial Analogue of a Theorem of R\'enyi}

\author{Kent E. Morrison}
\address{Mathematics Department, California Polytechnic State University, San Luis Obispo, CA 93407}
\email{kmorriso@calpoly.edu}

\subjclass[2000]{Primary 11T06; Secondary 11T55, 05A16}
\date{May 26, 2004}

\commby{XXXX}


\begin{abstract}
R\'enyi's result on the density of integers whose 
prime factorizations have excess multiplicity has an analogue for 
polynomials over a finite field.
\end{abstract}

\maketitle

Let $n=p_{1}^{\alpha_{1}}\cdots p_{r}^{\alpha_{r}}$ be the 
prime factorization of a positive integer $n$. Define the 
\textbf{excess} of $n$ to be $(\alpha_{1}-1) + \cdots + 
(\alpha_{r}-1)$, which is the difference between the total 
multiplicity $\alpha_{1}+\cdots + \alpha_{r}$ and the number of 
distinct primes in the factorization. An integer with excess 0 is 
also said to be \textbf{square-free}. Let $E_{k}$ denote the set of 
positive integers of excess $k$, $k=0,1,2,\ldots$. R\'enyi proved that 
the set $E_{k}$ has a density $d_{k}$ and that the sequence $\{ d_{k}\}$ 
has a generating function given by 
\[
   \sum_{k \geq 0}d_{k}z^{k}= \prod_{p} \left( 1-\frac{1}{p}\right)
                               \left( 1 + \frac{1}{p-z} \right) ,
\]
where the product extends over the primes. Recall that the density of 
a set of positive integers $E$ is the limit, if it exists,
\[
   \lim_{n\ra \infty} \frac{\#(E \cap \{1,2,\ldots,n \})}{n},
\]
which is the limiting probability that an integer from 1 to $n$ is in 
$E$.

The set of square-free integers is $E_{0}$ and setting $z=0$ in the 
generating function gives $d_{0}=\prod_{p}(1 - 1/p^{2})$, which is the 
well-known result that the density of square-free integers is 
$1/\zeta(2) = 6/\pi^{2}$. (This was first proved by Gegenbauer \cite{Gegenbauer85} in 1885. 
A clear, non-rigorous presentation is in \cite{Jones93}.)  By 
setting  $z=1$ one sees that $\sum_{k}d_{k}=1$, so that the density 
is countably additive on the specific partition of $\Z^{+}$ given by the $E_{k}$. 
R\'enyi's proof appeared in \cite{Renyi55}, but an alternative proof 
was given by Kac in \cite[pp. 64--71]{Kac59}. 

The aim of this paper is to derive an analogue of the generating 
function for polynomials in one variable over a finite field. Let 
$\fq$ be the field with $q$ elements and $\fq[x]$ the polynomial ring. 
The prime elements of $\fq[x]$ are the irreducible monic polynomials. 
Let $f$ be a monic polynomial with prime factorization 
$f=\pi_{1}^{\alpha_{1}}\cdots\pi_{r}^{\alpha_{r}}$, and define the 
excess of $f$ to be $(\alpha_{1}-1) + \cdots +(\alpha_{r}-1)$ just as 
for an integer. Let $e_{n,k}$ be the number of monic polynomials of 
degree $n$ and excess $k$. Define 
\[ d_{n,k} = \frac{e_{n,k}}{q^{n}},  \]
which is the probability that a monic polynomial of degree $n$ has 
excess $k$. Note that $d_{0,k}=0$ for $k > 0$. Then define the analogue 
of the density to be the  limiting ``probability" as the degree goes to infinity
\[   d_{k} = \lim_{n \ra \infty} d_{n,k}.  \]
Define $D(z) = \sum_{k \geq 0}d_{k}z^{k}$ to be the ordinary power 
series generating function of the sequence  $\{d_{k} \}$. Let $\norm f =q^{\deg f}$ be the norm of the polynomial $f$, which is the  cardinality of the residue ring $\fq[x]/(f)$. The main 
result of this paper is the following theorem concerning $D(z)$.
\begin{theorem} The generating function $D(z)$ has a factorization 
over the prime polynomials given by
\[
   D(z) = 
   \prod_{\pi}\left(1-\frac{1}{\norm \pi} \right)
               \left(1 + \frac{1}{\norm \pi-z} \right).
\]
\end{theorem}
\begin{proof}  We begin with the geometric series 
\[
   \frac{1}{1-qt} = \sum_{n \geq 0}q^{n}t^{n},
\]
which is the generating function for the number of monic polynomials 
of degree $n$. 
Then unique factorization in $\fq[x]$ allows us to factor the 
generating function formally
\begin{eqnarray}
   \frac{1}{1-qt} &=& \prod_{\pi}\sum_{j \geq 0}t^{j \deg \pi}  \label{prfact}\\
                  &=& \prod_{\pi} \frac{1}{1-t^{\deg \pi}}. 	 \nonumber	        
\end{eqnarray}
By grouping the primes of the same degree and letting 
$\nu_{i}$ denote the number of primes of degree $i$, we can rewrite 
the last line above as
\[  \frac{1}{1-qt} = \prod_{i \geq 1}\left(\frac{1}{1-t^{i}}\right)^{\nu_{i}}.
\]	
From this it follows that
\begin{equation} 
1-qt = 	\prod_{i \geq 1}\left({1-t^{i}}\right)^{\nu_{i} }   \label{prfact2} 
\end{equation}
as a formal power series. In the product on the right there is a finite number of terms for each power of $t$ so that the coefficients make sense. In fact, the coefficient of $t^n$ is 0 except for $n=0, 1$. However, considered as a function of a complex variable $t$, the product does not converge for all $t$. It does converges absolutely for $ |t| < 1/q$. This follows from consideration of the series $ \sum \nu_i t^i$ and the fact that $\nu_i$ is asymptotic to $q^i/i$. 

Next we define the two-variable generating function 
\[ E(t,z) = \sum_{n,k}e_{n,k}t^{n}z^{k}. 
\]
Modifying the factorization in (\ref{prfact}), we see that
\[
   E(t,z)=\prod_{\pi}(1+t^{\deg \pi} + t^{2 \deg \pi}z + \cdots + 
           t^{j \deg \pi}z^{j-1} + \cdots ).   
\]
Notice that the variable $z$ appears with a power that is equal to the
excess multiplicity.  That is, if
$f=\pi_{1}^{\alpha_{1}}\cdots\pi_{r}^{\alpha_{r}}$ then the product 
expansion of $E(t,z)$ has a term of the form 
$t^{\alpha_{1}\deg \pi_{1}}\cdots t^{\alpha_{r}\deg \pi_{r}}
 z^{\alpha_{1}-1}\cdots z^{\alpha_{r}-1}$. Sum the geometric 
 series in each factor to obtain the formal factorization
\[
   E(t,z)=\prod_{\pi}\left(1+ \frac{t^{\deg \pi}}{1-t^{\deg \pi}z} \right).
\]
Group the irreducibles by degree to get
\begin{equation}
  E(t,z)=\prod_{i}\left(1+ \frac{t^i}{1-t^i z} \right)^{\nu_i}. \label{prfact3}
\end{equation}
Now the product on the right converges absolutely if and only if the series 
\begin{equation}
 \sum_{i \geq 1} \nu_i \left|  \frac{t^i}{1-t^i z} \right|   \label{series}
 \end{equation}
converges. We claim that (\ref{series}) converges for $|t| < 1/q$ and $|z| < q$, 
because the denominators $|1-t^iz|$ are bounded away from 0 and $\nu_i$ is asymptotic to $q^i/i$. (Actually, it suffices that $\nu_i < q^i$.)

From (\ref{prfact2}) and (\ref{prfact3}) we get
\[ (1-qt) E(t,z) = \prod_{i \geq 1} \left( 1-t^{i}\right)^{\nu_i}
                          \prod_{i \geq 1}\left(1+ \frac{t^i}{1-t^i z} \right)^{\nu_i}  .\]
On the domain where both products converge absolutely we can combine the factors for each $i$ to get
\begin{equation} (1-qt) E(t,z) = \prod_{i \geq 1} \left( 1-t^{i}\right)^{\nu_i}
                          \left(1+ \frac{t^i}{1-t^i z} \right)^{\nu_i} .  \label{product}
\end{equation}
By multiplying the factors together we can see that the absolute convergence of the infinite product depends on the convergence of the series
\[  \sum_i \nu_i \left|  \frac{t^{2i}z - t^{2i}}{1-t^i z}   \right|  .     
\]
Then reasoning along the same lines as before we see that this series converges for $|t^2| < q $ and $|z| < \sqrt{q}$. In particular, the product converges for $t=1/q$, and so after carrying out the multiplication of the left side of (\ref{product}) we arrive at
\[
   \sum_{n,k}(e_{n,k}- q e_{n-1,k})t^{n}z^k =  \prod_{i \geq 1} \left( 1-t^{i}\right)^{\nu_i}
                          \left(1+ \frac{t^i}{1-t^i z} \right)^{\nu_i} .
\]                           
We evaluate this at $t=1/q$ to get
\[
    \sum_{n,k}(e_{n,k}- q e_{n-1,k})(1/q)^{n}z^k =  \prod_{i \geq 1} \left( 1-(1/q)^{i}\right)^{\nu_i}
                          \left(1+ \frac{(1/q)^i}{1-(1/q)^i z} \right)^{\nu_i} .
\]
The coefficient of $z^k$ is the sum $\sum_{n \geq 1} ( e_{n,k}/q^n - e_{n-1,k}/q^{n-1})$. This  telescopes to give
\[ \lim_{n \ra \infty} \frac{ e_{n,k}}{q^n } = \lim_{n \ra \infty} d_{n,k},\] which is the definition of $d_k$, and so we  have
\[   D(z)=\sum_{k}d_{k}z^{k}  =  \prod_{i \geq 1} \left( 1-(1/q)^{i}\right)^{\nu_i}
                          \left(1+ \frac{(1/q)^i}{1-(1/q)^i z} \right)^{\nu_i}  .
 \]                        
 Finally, we write the product by indexing over the prime polynomials $\pi$ and note that the norm of $\pi$ is $\norm \pi = q^{\deg \pi}$. With this we have the generating function for the $d_k$ in the form that is most directly analogous to R\'enyi's theorem.
 \begin{eqnarray*}
   D(z)&=&\prod_{\pi}\left(1-(1/q)^{\deg \pi} \right)
       \left(1+ \frac{(1/q)^{\deg \pi}}{1-(1/q)^{\deg \pi}z} 
       \right)  \nonumber \\
       &=& \prod_{\pi}\left(1-\frac{1}{\norm \pi} \right)
               \left(1 + \frac{1}{\norm \pi-z} \right). 
               \label{dkgenfct}
\end{eqnarray*}
\end{proof}
The coefficient $d_0$ is the limiting ``probability'' that a monic polynomial is square-free.   To develop the analogy with the density of the square-free integers given by $d_0$ in R\'enyi's generating function, we use the zeta function of $\fq[x]$ (i.e. the zeta function of the affine line over $\fq$) 
\[ \zeta(s) = \frac{1}{1-q^{-s}}  ,\]
which immediately comes from the definition 
\[ \zeta(s) = \sum_{\mathfrak{a}} \frac{1}{(\norm \mathfrak{a})^s}, \]
where the sum is over all ideals of $\fq[x]$ and the norm $\norm \mathfrak{a}$ is the cardinality of the residue ring $\fq[x]/\mathfrak{a}$.  It has a factorization over the prime ideals (i.e. irreducible polynomials)
\begin{eqnarray*}
\zeta(s) &=& \prod_{\pi} \frac{1}{1-(\norm \pi)^{-s}}  \\
              &=& \prod_{i \geq 1}\left( \frac{1}{1-q^{-is}} \right)^{\nu_i}.
\end{eqnarray*}
\begin{corollary}
   $d_{0}= \frac{1}{\zeta(2)}=1-\frac{1}{q}.$
\end{corollary}
\begin{proof}
We have 
\[ 
   d_0 = D(0) = \prod_{i \geq 1} \left( 1- \frac{1}{q^{2i}}  \right)^{\nu_i}.
\]
Then in (\ref{prfact2}) we may let $t=1/q^2$, because the product converges for $|t| < 1/q$, to obtain
\[
   1-\frac{1}{q}= \prod_{i \geq 1} \left(1-\frac{1}{q^{2i}}\right)^{\nu_i}.
\]
Notice that the product is $1/\zeta(2)$.
\end{proof}
Corollary 1 can be obtained as a special case of much more general results on square-free values of polynomials in one or more variables from the work of Ramsay \cite{Ramsay92} and Poonen \cite{Poonen03}.
It turns out that for $n \geq 2$, the value of $d_{n,0}$ is $1-1/q$. This can be seen by finding the coefficients $e_{n,0}$, which count the number of monic, square-free polynomials of degree $n$. These polynomials can be counted directly; see, for example, \cite{Carlitz56}.
\begin{corollary}
   The number of square-free monic polynomials of degree $n \geq 2$ is 
   $q^{n}-q^{n-1}$.
\end{corollary}
\begin{proof}  The generating function $\sum_{n \geq 0} e_{n,0}t^n = E(t,0)$. From  (\ref{prfact3}) we see that
\[
  E(t,0)=\prod_{i \geq 1} (1 + t^i)^{\nu_i}.
  \]
Using (\ref{prfact2}) we see that
\begin{eqnarray*}
   E(t,0)(1-qt) &=&\prod_{i \geq 1}(1 + t^i)^{\nu_i}(1- t^i)^{\nu_i}    \\
      &=&  \prod_{i \geq 1}(1 - t^{2i})^{\nu_i} \\
      &=& 1-qt^2.
\end{eqnarray*}
Therefore,
\[  E(t,0) = \frac{1-qt^2}{1-qt}  , \]
from which it follows that $e_{n,0} = q^n -q^{n-1}$ for $n \geq 2$.
 \end{proof}
\vspace{2mm}
From the expression 
\[ D(z)=\prod_{i \geq 1}\left(1-\frac{1}{q^{i}} 
\right)^{\nu_{i}}
               \left(1 + \frac{1}{q^{i}-z} \right)^{\nu_{i}}
\]
we can see that $D(z)$ has poles at $z=q^{i}$ of multiplicity $\nu_{i}$. 
In particular the pole at $z=q$ has multiplicity $q-1$. Elementary 
analysis of the singularity there, along the lines of Kac \cite{Kac59} in his discussion of R\'enyi's result, enables us to describe the 
asymptotic behavior of the $d_{k}$ as $k \ra \infty$.
\begin{corollary}
   As $k$ goes to infinity, $d_{k}$ is asymptotic to
   \[ A\frac{k^{q-2}}{q^{k}} , \]
where the constant $A$ is given by   
   \[ A= \frac{1}{(q-2)!}\left(\frac{1}{q}-\frac{1}{q^{2}} \right)^{q-1}
   \prod_{i \geq 2}\left(1-\frac{1}{q^{i}}\right)^{\nu_{i}}
                   \left(1-\frac{1}{q^{i}-q}\right)^{\nu_{i}}.
   \]
\end{corollary}

One may contrast this asymptotic result with the classical case of  R\'enyi. Although the generating functions have clearly analogous form, the generating function for the number-theoretic version has only a simple pole $z=2$, which is the pole of smallest absolute value. The asymptotic analysis shows that
 \[d_{k} \sim 
\frac{\delta}{2^{k}} , \]
where \[ \delta = \frac{1}{4}\prod_{p \geq 
3}\frac{(p-1)^{2}}{p(p-2)}. \]

The referee has observed that Theorem 1 of this article can be extended naturally to function fields over finite fields by using $S$-zeta functions and their residues at $t=1/q$.  Let $K$ be a function field over the constant field $\fq$. Let $S$ be a finite, non-empty set of places on $K$ and let $\mathcal{O}_{K,S}$ denote the ring of $S$-integers of $K$. Then for every integer $k \geq 0$, the density $d_{k,S}$ of ideals in $\mathcal{O}_{K,S}$ with excess $k$ exists and the following analytic identity holds:
\[ \sum_{k \geq 0} d_kz^k = \prod_{v \notin S}\left( 1 - \frac{1}{\norm v} \right)\left( 1+\frac{1}{\norm v - z}  \right)  ,  \]
where for each place $v$ on $K$ the norm $\norm v = q^{\deg v}$ is the cardinality of the residue field at $v$.

Finally, the referee has pointed out that by generalizing Kac's proof of R\'enyi's theorem
 \cite[pp. 64--71]{Kac59}, there should also be an analogue of the theorem for the density of ideals with excess $k$ in the ring of algebraic integers (or ring of $S$-integers) of any number field.

\end{document}